\documentclass[11pt]{article}

\usepackage{latexsym}
\usepackage{graphicx}
\usepackage{float}
\usepackage{amsmath}
\usepackage{amssymb}
\usepackage{lscape}
\usepackage{mathtools}
\usepackage{amsthm}
\usepackage{tikz}    
\usepackage{mathptmx}

\begin{document}
\bibliographystyle{plain}
\floatplacement{table}{H}
\newtheorem{definition}{Definition}[section]
\newtheorem{lemma}{Lemma}[section]
\newtheorem{theorem}{Theorem}[section]
\newtheorem{corollary}{Corollary}[section]
\newtheorem{proposition}{Proposition}[section]

\newcommand{\sni}{\sum_{i=1}^{n}}
\newcommand{\snj}{\sum_{j=1}^{n}}
\newcommand{\smj}{\sum_{j=1}^{m}}
\newcommand{\sumjm}{\sum_{j=1}^{m}}
\newcommand{\bdis}{\begin{displaymath}}
\newcommand{\edis}{\end{displaymath}}
\newcommand{\beq}{\begin{equation}}
\newcommand{\eeq}{\end{equation}}
\newcommand{\beqn}{\begin{eqnarray}}
\newcommand{\eeqn}{\end{eqnarray}}
\newcommand{\simleq}{\stackrel{<}{\sim}}
\newcommand{\sep}{\;\;\;\;\;\; ; \;\;\;\;\;\;}
\newcommand{\real}{\mbox{$ I \hskip -4.0pt R $}}
\newcommand{\complex}{\mbox{$ I \hskip -6.8pt C $}}
\newcommand{\integ}{\mbox{$ Z $}}
\newcommand{\realn}{\real ^{n}}
\newcommand{\sqrn}{\sqrt{n}}
\newcommand{\sqrtwo}{\sqrt{2}}
\newcommand{\prf}{{\bf Proof. }}

\newcommand{\onehlf}{\frac{1}{2}}
\newcommand{\donehlf}{\dfrac{1}{2}}
\newcommand{\thrhlf}{\frac{3}{2}}
\newcommand{\fivhlf}{\frac{5}{2}}
\newcommand{\onethd}{\frac{1}{3}}
\newcommand{\lb}{\left ( }
\newcommand{\lcb}{\left \{ }
\newcommand{\lsb}{\left [ }
\newcommand{\labs}{\left | }
\newcommand{\rb}{\right ) }
\newcommand{\rcb}{\right \} }
\newcommand{\rsb}{\right ] }
\newcommand{\rabs}{\right | }
\newcommand{\lnm}{\left \| }
\newcommand{\rnm}{\right \| }
\newcommand{\lambdab}{\bar{\lambda}}
%
%
\newcommand{\xj}{x_{j}}
\newcommand{\xjb}{\bar{x}_{j}}

\newcommand{\anmell}{a_{n-\ell}}
\newcommand{\anmellmo}{a_{n-\ell-1}}
\newcommand{\anmkmo}{a_{n-k-1}}
\newcommand{\anmkmell}{a_{n-k-\ell}}
\newcommand{\anmkmellmo}{a_{n-k-\ell-1}}
\newcommand{\azero}{a_{0}}
\newcommand{\aone}{a_{1}}
\newcommand{\atwo}{a_{2}}
\newcommand{\ath}{a_{3}}
\newcommand{\afr}{a_{4}}
\newcommand{\afv}{a_{5}}
\newcommand{\asx}{a_{6}}
\newcommand{\anmo}{a_{n-1}}
\newcommand{\anmt}{a_{n-2}}
\newcommand{\an}{a_{n}}
\newcommand{\aj}{a_{j}}
\newcommand{\ajpo}{a_{j+1}}
\newcommand{\ajpt}{a_{j+2}}
\newcommand{\ajmo}{a_{j-1}}
\newcommand{\ajmt}{a_{j-2}}
\newcommand{\ak}{a_{k}}
\newcommand{\akpo}{a_{k+1}}
\newcommand{\akpt}{a_{k+2}}
\newcommand{\akmo}{a_{k-1}}
\newcommand{\akmt}{a_{k-2}}
\newcommand{\aell}{a_{\ell}}

\newcommand{\bone}{b_{1}}
\newcommand{\btwo}{b_{2}}
\newcommand{\defeq}{\coloneqq}

\newcommand{\zj}{z^{j}}          
\newcommand{\zk}{z^{k}}          
\newcommand{\zkpm}{z^{k+m}}          
\newcommand{\zm}{z^{m}}          

\newcommand{\matrixspace}{\;\;}
\newcommand{\ellone}{\ell_{1}} \newcommand{\elltwo}{\ell_{2}}

\newcommand{\ronenk}{r_{1}}
\newcommand{\rtwonk}{r_{2}}
\newcommand{\obar}{\widebar{O}}
\newcommand{\Lbar}{\widebar{L}}
\newcommand{\dltaone}{\delta_{1}}
\newcommand{\dltatwo}{\delta_{2}}
\newcommand{\dltatld}{\tilde{\delta}}
\newcommand{\tauone}{\tau_{1}}
\newcommand{\tautwo}{\tau_{2}}
\newcommand{\muone}{\mu_{1}}
\newcommand{\mutwo}{\mu_{2}}
\newcommand{\xb}{\bar{x}}
\newcommand{\qone}{q_{1}}
\newcommand{\qtwo}{q_{2}}
\newcommand{\ub}{\bar{u}}
\newcommand{\cmm}{\complex^{m \times m}}
\newcommand{\opdsk}{\mathcal{O}}
\newcommand{\cldsk}{\widebar{\mathcal{O}}}

\newcommand{\sqrtaone}{\sqrt{\a1}}
\newcommand{\sqrtatwo}{\sqrt{\atwo}}

\newcommand{\rhoone}{\rho_{1}}
\newcommand{\rhotwo}{\rho_{2}}
\newcommand{\xone}{x_{1}}
\newcommand{\xtwo}{x_{2}}
\newcommand{\xthr}{x_{3}}
\newcommand{\xfor}{x_{4}}
\newcommand{\xfiv}{x_{5}}

\newcommand{\lnorm}{\left \|}
\newcommand{\rnorm}{\right \|}
\newcommand{\lnrm}{\biggl | \biggl |}
\newcommand{\rnrm}{\biggr |\biggr |}

\newcommand{\recipp}{p^{\protect \#}}
\newcommand{\recipP}{P^{\protect \#}}

\newcommand{\absatwo}{|\atwo|}
\newcommand{\absqrtatwo}{\sqrt{|\atwo|}}

\newcommand{\xstar}{x^{\ast}}
\newcommand{\xkova}{\hat{x}}

\newcommand{\alfazero}{\alpha_{0}}
\newcommand{\alfaone}{\alpha_{1}}
\newcommand{\alfatwo}{\alpha_{2}}
\newcommand{\alfathr}{\alpha_{3}}
\newcommand{\alfath}{\alpha_{3}}
\newcommand{\alfafr}{\alpha_{4}}
\newcommand{\alfafv}{\alpha_{5}}
\newcommand{\alfasx}{\alpha_{6}}
\newcommand{\alfanmo}{\alpha_{n-1}}
\newcommand{\alfanmt}{\alpha_{n-2}}
\newcommand{\alfan}{\alpha_{n}}
\newcommand{\alfaj}{\alpha_{j}}
\newcommand{\alfai}{\alpha_{i}}
\newcommand{\alfak}{\alpha_{k}}
\newcommand{\alfajpo}{\alpha_{j+1}}
\newcommand{\alfakpo}{\alpha_{k+1}}
\newcommand{\betaone}{\beta_{1}}
\newcommand{\betatwo}{\beta_{2}}

\newcommand{\Azero}{A_{0}}
\newcommand{\Aone}{A_{1}}
\newcommand{\Atwo}{A_{2}}
\newcommand{\Ath}{A_{3}}
\newcommand{\Afr}{A_{4}}
\newcommand{\Afv}{A_{5}}
\newcommand{\Asx}{A_{6}}
\newcommand{\Anmo}{A_{n-1}}
\newcommand{\Anmt}{A_{n-2}}
\newcommand{\An}{A_{n}}
\newcommand{\Aj}{A_{j}}
\newcommand{\Aell}{A_{\ell}}
\newcommand{\Aellpo}{A_{\ell+1}}
\newcommand{\Aellmo}{A_{\ell-1}}
\newcommand{\Aellmk}{A_{\ell-k}}
\newcommand{\Aellpk}{A_{\ell+k}}
\newcommand{\Ajmk}{A_{j-k}}

\newcommand{\lnq}{<}
\newcommand{\gnq}{>}

\newcommand{\abels}{(a,b) \in S}

\newcommand{\eneskakk}{Enestr\"{o}m-Kakeya}
\newcommand{\realpart}{\operatorname{Re}}
\newcommand{\imagpart}{\operatorname{Im}}
\newcommand{\bmult}{\begin{multline}}
\newcommand{\emult}{\end{multline}}
\newcommand{\phij}{\varphi_{j}}

\newcommand{\uazerou}{u^{*}A_{0}u}
\newcommand{\uaoneu}{u^{*}A_{1}u}

\newcommand{\AH}{A_{H}}
\newcommand{\AS}{A_{S}}
\newcommand{\leftangle}{\left \langle}
\newcommand{\rightangle}{\right \rangle}
\newcommand{\ustar}{u^{*}}
\newcommand{\uAu}{\ustar A u}
\newcommand{\uAHu}{\ustar \AH u}
\newcommand{\uASu}{\ustar \AS u}
\newcommand{\uMu}{\ustar M u}
\newcommand{\lambdamin}{\lambda_{min}}
\newcommand{\lambdamax}{\lambda_{max}}

\newcommand{\Th}{T_{H}}
\newcommand{\Ts}{T_{S}}
\newcommand{\uTu}{\leftangle Tu , u \rightangle}
\newcommand{\uQu}{\leftangle Qu , u \rightangle}
\newcommand{\uthu}{\leftangle \Th u , u \rightangle}
\newcommand{\utsu}{\leftangle \Ts u , u \rightangle}

\newcommand{\utoneu}{\leftangle T_{1}u,u \rightangle}
\newcommand{\uttwou}{\leftangle T_{2}u,u \rightangle}

\newcommand{\bh}{\mathcal{B}(H)}
\newcommand{\wtclos}{\overline{W(T)}}

%
%
%
%

\begin{center}
\large
{\bf AN OCTAGON CONTAINING THE NUMERICAL RANGE OF A BOUNDED LINEAR OPERATOR}
\vskip 0.5cm
\normalsize
A. Melman \\
Department of Applied Mathematics \\
School of Engineering, Santa Clara University  \\
Santa Clara, CA 95053  \\
e-mail : amelman@scu.edu \\
\vskip 0.5cm
\end{center}

\begin{abstract}
A polygon is derived that contains the numerical range of a bounded linear operator on a complex Hilbert space, using only norms.
In its most
general form, the polygon is an octagon, symmetric with respect to the origin, and tangent to the closure of the numerical range in at least four points
when the spectral norm is used. 
\vskip 0.15cm
{\bf AMS(MOS) subject classification :} 47A12, 47L30, 15A60, 65H17
\end{abstract}

%
%
%
%
\section{Introduction}
\label{introduction}

The \emph{numerical range} of $T \in \mathcal{B}(H)$, the algebra of bounded linear operators on a complex Hilbert space~$H$, equipped with the inner product 
$\leftangle .,. \rightangle$, is the subset of~$\mathbb{C}$, defined by
\bdis
W(T) = \lcb \uTu : u \in H \, , \, \|u\|=1 \rcb \; , 
\edis
where $\|u\|^{2} = \leftangle u,u \rightangle$.
Also referred to as the \emph{field of values}, it plays an important role in several fields of mathematics and engineering.
By the Toeplitz-Hausdorff theorem, $W(T)$ is a convex set. A related quantity, the \emph{numerical radius}, is defined as
$ w(T)=\sup_{\|u\|=1} |\uTu| $. 

The numerical range can be enclosed by a polygonal envelope (for matrices, but easily generalized to bounded operators, see\cite[Section 1.5]{HJ_TIMA}), 
although this requires the computation of eigenvalues and corresponding
eigenvectors, which is impractical when matrix sizes are large or in cases where the matrix is only implicitly defined. 

Our purpose is to
enclose the numerical range in an easily computable region (in its most general form an octagon) using only norms and avoiding the computation 
of spectral or spectral-related quantitities. On the one hand, this leads to a cruder approximation than could be obtained by using                 
spectral information, but on the other, it is faster and much simpler. It will depend on the application whether accuracy or computational simplicity
is preferable, but such matters are beyond our scope here. 

To begin, we briefly review a few basic properties of $\mathcal{B}(H)$ and the numerical range, as can be found in any standard text on these subjects
(e.g., \cite{AG}, \cite{GR}, \cite{Weidmann}).
We denote by $T^{*}$ the adjoint of $T \in \bh$, defined by $\leftangle Tu,u\rightangle = \leftangle u,T^{*}u \rightangle$, $u \in H$. 
An operator~$T$ is self-adjoint if $T=T^{*}$. The Cartesian decomposition of $T \in \bh$ is given by $T=\Ts + i\Ts$,
where $\Th$ and $\Ts$ are self-adjoint bounded operators defined as 
\beq 
\nonumber         
\Th = \dfrac{1}{2} \lb T + T^{*} \rb \;\; \text{and} \;\; \Ts = \dfrac{1}{2i} \lb T - T^{*} \rb \; .  
\eeq   
It follows from this decomposition that $\uTu = \uthu + i \utsu$, $u \in H$, with $\uthu , \utsu \in \mathbb{R}$.  

The spectral norm of~$T \in \bh$ is defined as $\|T\|_{\sigma} = \sup_{\|u\|=1} \|Tu\|$. 
There exist several upper bounds for~$w(T)$, expressed in terms of the spectral norm: first, as an immediate consequence of the 
definition of~$w(T)$, one has the standard bound 
\beq
\label{bnd_classical}
w(T) \leq \|T \|_{\sigma}  \; .
\eeq
However, this bound is not necessarily satisfied when the norm is different from the spectral norm:
a finite dimensional counterexample of a $2 \times 2$ matrix (a bounded linear operator on $\mathbb{C}^{2}$) with the matrix 1-norm is given by
\bdis
\dfrac{3+\sqrt{3}}{4} =
\left |
\begin{pmatrix}
\sqrt{3}/2 \\
1/2 \\
\end{pmatrix}^{*}
\begin{pmatrix}
1 & 1 \\
0   & 0   \\
\end{pmatrix}
\begin{pmatrix}
\sqrt{3}/2 \\
1/2 \\
\end{pmatrix}
\right |
>               
\left \|
\begin{pmatrix}
1 & 1 \\
0   & 0   \\
\end{pmatrix}
\right \|_{1}
= 1
\; .
\edis
Two recent improvements of the bound in~(\ref{bnd_classical}) are the following:
\begin{eqnarray}
& & w(T) \leq \dfrac{1}{2} \lb \|T\|_{\sigma} + \|T^{2}\|_{\sigma}^{1/2} \rb \;\; \text{from~\cite{Kit1},}  \label{kittaneh_ineq1} \\
& & w(T) \leq \lb \dfrac{\|T^{*}T+TT^{*}\|_{\sigma}}{2} \rb ^{1/2} \;\; \text{from~\cite{Kit2}.}  \label{kittaneh_ineq2} 
\end{eqnarray}

When $T \in \bh$ is self-adjoint, then $w(T) = \|T\|_{\sigma}$, $w(T) \leq \|T \|$ for any norm, 
and $\overline{W(T)} \subseteq [-\|T\|_{\sigma},\|T\|_{\sigma}]$, where $\wtclos$ is the closure of~$W(T)$.
Throughout, we denote the real and imaginary parts of a complex number~$z$ by $\Re z$ and $\Im z$, respectively.

We now derive the enclosing octagon mentioned earlier.

%
%
%
%
\section{An octagon containing the numerical range}
\label{octagon}   

\begin{figure}[H]
\begin{center}
\begin{tikzpicture}[scale=0.25]
        \tikzstyle{every path}=[draw, line width=0.5mm, color=black];
        \coordinate (v1) at (-10,8) {};
        \coordinate (v2) at (10,8) {};
        \coordinate (v3) at (10,-8) {};
        \coordinate (v4) at (-10,-8) {};
        \coordinate (p1) at (6,8) {};
        \coordinate (p2) at (10,4) {};
        \coordinate (p3) at (10,-2) {};
        \coordinate (p4) at (4,-8) {};
        \coordinate (p5) at (-6,-8) {};
        \coordinate (p6) at (-10,-4) {};
        \coordinate (p7) at (-10,2) {};
        \coordinate (p8) at (-4,8) {};


        \draw[opacity=0.5,fill=gray] (v1.center)--(v2.center);
        \draw[opacity=0.5,fill=gray] (v2.center)--(v3.center);
        \draw[opacity=0.5,fill=gray] (v3.center)--(v4.center);
        \draw[opacity=0.5,fill=gray] (v4.center)--(v1.center);

        \draw[opacity=0.5,fill=gray] (p1.center)--(p2.center);
        \draw[opacity=0.5,fill=gray] (p2.center)--(p3.center);
        \draw[opacity=0.5,fill=gray] (p3.center)--(p4.center);
        \draw[opacity=0.5,fill=gray] (p4.center)--(p5.center);
        \draw[opacity=0.5,fill=gray] (p5.center)--(p6.center);
        \draw[opacity=0.5,fill=gray] (p6.center)--(p7.center);
        \draw[opacity=0.5,fill=gray] (p7.center)--(p8.center);
        \draw[opacity=0.5,fill=gray] (p8.center)--(p1.center);
        \fill [opacity=0.45,gray] (p1) \foreach \i in {1,...,8}{ -- (p\i) } -- cycle;     
        
        \draw[dashed,very thick] (0.5,13.5) -- (13.5,0.5);
        \draw[dashed,very thick] (-0.5,-13.5) -- (-13.5,-0.5);
        \draw[dashed,very thick] (-1.5,-13.5) -- (13.5,1.5);
        \draw[dashed,very thick] (-1.5,-13.5) -- (13.5,1.5);
        \draw[dashed,very thick] (-13.5,-1.5) -- (1.5,13.5);
        \foreach \i in{1,...,8} \fill [very thick] (p\i) circle [radius=.25]; 
        \draw[solid,very thick] (0,-15) -- (0,15);
        \draw[solid,very thick] (-15,0) -- (15,0);
        \node at (5.25,11) {{\Large $\mathbf L_{1}$}};
        \node at (-3.25,11) {{\Large $\mathbf L_{4}$}};
        \node at (-5,-11) {{\Large $\mathbf L_{2}$}};
        \node at (3.5,-11) {{\Large $\mathbf L_{3}$}};
        \fill [very thick] (10,8) circle [radius=.05]; 
        \draw [thin] (10,8) circle [radius=.25]; 
        \node at (15.5,8) {{\large ${\mathbf (\|\Th\|,\|\Ts\|)}$}};
        \fill [very thick] (10,-8) circle [radius=.05]; 
        \draw [thin] (10,-8) circle [radius=.25]; 
        \node at (16,-8) {{\large ${\mathbf (\|\Th\|,-\|\Ts\|)}$}};
        \fill [very thick] (-10,-8) circle [radius=.05]; 
        \draw [thin] (-10,-8) circle [radius=.25]; 
        \node at (-17.15,-8) {{\large ${\mathbf (-\|\Th\|,-\|\Ts\|)}$}};
        \fill [very thick] (-10,8) circle [radius=.05]; 
        \draw [thin] (-10,8) circle [radius=.25]; 
        \node at (-16.5,8) {{\large ${\mathbf (-\|\Th\|,\|\Ts\|)}$}};
        \fill [very thick] (0,0) circle [radius=.05]; 
        \draw [thin] (0,0) circle [radius=.25]; 
        \node at (-1,-1) {{\Large ${\mathbf 0 }$}};
\end{tikzpicture}
\caption{$\mathcal{R}(T)$ and $\mathcal{P}_{\alpha\beta\gamma\delta}$ for Theorem~\ref{theorem_rectangles}.}
\label{fig1}                             
\end{center}
\end{figure}
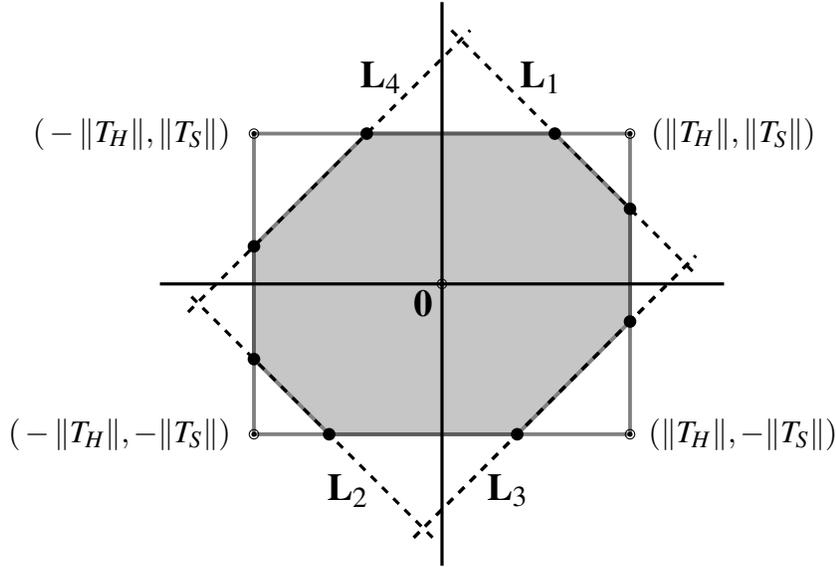

The following theorem forms the basis for the construction of a polygon containing the numerical range.
%
%
%
%
\begin{theorem}     
\label{theorem_rectangles}
Let $T \in \mathcal{B}(H)$ have the Cartesian decomposition $T=\Th + i\Ts$, let $T \neq c\, Q$ 
for any $c \in \mathbb{C}$ and $Q \in \mathcal{B}(H)$ with $Q=Q^{*}$,
and let $\alpha,\beta,\gamma,\delta > 0$. For any norm, define 
the rectangle~$\mathcal{R}(T)$ and the parallelogram~$\mathcal{P}_{\alpha\beta\gamma\delta}(T)$, both centered at the origin in the complex plane, by
\beq
\nonumber      
\mathcal{R}(T) = \Bigl \{ x + iy : x,y \in \mathbb{R}\, , \, |x| \leq \|\Th\| \; \text{and} \;\, |y| \leq \|\Ts\| \Bigr \} 
\eeq
and
\begin{multline}
\label{Wlines}
\mathcal{P}_{\alpha\beta\gamma\delta}(T) = \Bigl \{ x + iy : x,y \in \mathbb{R} \, , \, |\alpha x + \beta y| \leq \| \alpha \Th + \beta \Ts \| \\
\; \text{and} \;\, |\gamma x - \delta y| \leq \| \gamma \Th - \delta \Ts \| \Bigr \}  . 
\end{multline}
Then the following holds.
\begin{itemize}
\item[{\bf (1)}]
$\wtclos \subseteq \mathcal{R}(T) \cap \mathcal{P}_{\alpha\beta\gamma\delta}(T)$.
\item[{\bf (2)}]
The corner points of the rectangle $\mathcal{R}(T)$ either lie
outside \newline $\mathcal{R}(T) \cap \mathcal{P}_{\alpha\beta\gamma\delta}(T)$ or on the boundary of this intersection.
\item[{\bf (3)}]
If the spectral norm is used to construct~$\mathcal{R}(T)$ and~$\mathcal{P}_{\alpha\beta\gamma\delta}(T)$, then each side 
or its opposing side of $\mathcal{R}(T)$ and $\mathcal{P}_{\alpha\beta\gamma\delta}(T)$ is tangent to~$\wtclos$, where the disjunction
is inclusive.
\end{itemize}
If $T=c\, Q$ for $c \in \mathbb{C}$ and $Q \in \mathcal{B}(H)$ with $Q=Q^{*}$,
then~$\wtclos$ is contained in the closed line segment determined by the endpoints $\pm c \, \|Q\|_{\sigma}$, at least one of which is 
a boundary point of~$\wtclos$ if the norm is the spectral norm. 
\end{theorem}
\begin{proof}
Consider $T \in \bh$ that is not a complex multiple of a self-adjoint operator.
Since $\uTu = \uthu + i\utsu$ and $\Th$ and $\Ts$ are self-adjoint, we have that 
\bdis
\left | \Re \uTu \right | = |\uthu|  \leq \|\Th\| \;\; \text{and} \;\; \left | \Im \uTu \right | = |\utsu|\leq \|\Ts\| \; ,
\edis
so that $\uTu \in \mathcal{R}(T)$, and, since $\mathcal{R}(T)$ is closed, the limit points of any sequence $\lcb \leftangle T u_{n}, u_{n} \rightangle \rcb$,
$\|u_{n}\|=1$, 
also lie in~$\mathcal{R}(T)$.
Moreover, for any $\alpha, \beta \in \mathbb{R}$, 
\begin{multline}
\nonumber
\alpha \uthu + \beta \utsu = \leftangle \lb \alpha \Th + \beta \Ts \rb u,u \rightangle \\ 
\; \Longrightarrow \;
|\alpha \uthu + \beta \utsu| \leq \| \alpha \Th + \beta \Ts \| \; ,
\end{multline}
which is equivalent to 
\bdis
|\alpha \Re\uTu + \beta \Im\uTu| \leq \| \alpha \Th + \beta \Ts \| \; .
\edis
The second inequality in~(\ref{Wlines}) follows analogously for any $\gamma,\delta \in \mathbb{R}$.
When $\alpha, \beta, \gamma,\delta >0$, then the inequalities in~(\ref{Wlines}) define the closed parallelogram 
$\mathcal{P}_{\alpha\beta\gamma\delta}(T)$ centered at the origin, and bounded by the lines $L_{j}$ ($j=1,2,3,4$), defined, 
after the usual identification of $\mathbb{C}$ with $\mathbb{R}^{2}$, by
\begin{eqnarray*}
& & L_{1} : \alpha x + \beta y = \|\alpha \Th + \beta \Ts \| 
\;\; , \;\; 
L_{2} (x,y) : \alpha x + \beta y = -\|\alpha \Th + \beta \Ts \| \; , \\ 
& & L_{3} : \gamma x - \delta y = \|\gamma \Th - \delta \Ts \|
\;\; , \;\; 
L_{4} (x,y) : \gamma x - \delta y = -\|\gamma \Th - \delta \Ts \| \; , 
\end{eqnarray*}
as illustrated in Figure~\ref{fig1}. The lines $L_{j}$ define a nondegenerate parallelogram because their right-hand sides never vanish, as the latter would
imply that $\Th$ and $\Ts$ are multiples of each other, and this was explicitly excluded by the condition that $T \neq c \, Q$ for a self-adjoint
operator~$Q$. This means that $\uTu$ and the limit points of any sequence $\lcb \leftangle T u_{n}, u_{n} \rightangle \rcb$,
$\|u_{n}\|=1$, 
lie in $\mathcal{P}_{\alpha\beta\gamma\delta}(T)$ as well and the first part of the theorem follows.

We prove the second part for the upper and lower right-hand corners of~$\mathcal{R}(T)$ as the result then follows for the remaining corner points 
from the symmetry with respect to the origin of both $\mathcal{R}(T)$ and~$\mathcal{P}_{\alpha\beta\gamma\delta}(T)$.
For the upper right-hand corner $\lb \|\Th \| , \|\Ts \| \rb$, we obtain with $L_{1}$:
\begin{multline}
\nonumber 
\alpha \|\Th\| + \beta   \|\Ts \| \geq \|\alpha \Th + \beta \Ts \| \\ 
\Longrightarrow \lb \|\Th \| , \|\Ts \| \rb \in \partial\mathcal{P}_{\alpha\beta\gamma\delta}(T) \;\; \text{OR} \;
\notin \mathcal{P}_{\alpha\beta\gamma\delta}(T) \; ,   
\end{multline}   
whereas for the lower right-hand corner $\lb \|\Th \| , -\|\Ts \| \rb$, we obtain with $L_{3}$:
\begin{multline}
\nonumber 
\gamma \|\Th\|  + \delta \|\Ts \| \geq \|\gamma \Th - \delta \Ts \|  \\
\Longrightarrow \lb \|\Th\| , -\|\Ts\| \rb \in \partial\mathcal{P}_{\alpha\beta\gamma\delta}(T) \;\; \text{OR} \;
\notin \mathcal{P}_{\alpha\beta\gamma\delta}(T) \; ,    
\end{multline}    
and the second part of the proof follows.

For the last part of the proof, where the norm is assumed to be the spectral norm, we first consider the self-adjoint operator $\Th$, which satisfies 
$\|\Th\|_{\sigma}=\sup_{\|u\|=1} |\uthu|$. From this it follows that there exists a sequence $\lcb u_{n} \rcb$ in~$H$ with
$\|u_{n}\|=1$, such that 
\bdis
\|\Th\|_{\sigma}=\lim_{n \rightarrow \infty} |\leftangle \Th u_{n},u_{n} \rightangle |\; .
\edis
Therefore, the real sequence $\lcb \leftangle \Th u_{n},u_{n} \rightangle\rcb$ contains a subsequence 
$\lcb \leftangle \Th v_{n},v_{n} \rightangle\rcb$ that converges either to $\|\Th\|_{\sigma}$ or $-\|\Th\|_{\sigma}$.
Since $\lcb \leftangle \Th v_{n},v_{n} \rightangle\rcb = \lcb \Re \leftangle T v_{n},v_{n} \rightangle\rcb$, this means that
$\lcb \leftangle T v_{n},v_{n} \rightangle\rcb$ converges to the left or right side of~$\mathcal{R}(T)$, which is then  
necessarily tangent to~$\wtclos$. An analogous argument for the self-adjoint operator~$\Ts$ and a convergent sequence 
$\lcb \Im \leftangle T r_{n},r_{n} \rightangle\rcb$ shows that the top or bottom side of~$\mathcal{R}(T)$ is tangent to~$\wtclos$.
In the case of the self-adjoint operator~$\alpha\Th + \beta\Ts$, one similarly obtains with the help of a sequence 
$\lcb \alpha \Re \leftangle T s_{n},s_{n} \rightangle + \beta \Im \leftangle T s_{n},s_{n} \rightangle \rcb$ 
that the top right or bottom left side of~$\mathcal{P}_{\alpha\beta\gamma\delta}(T)$ is tangent to~$\wtclos$, 
and an analogous argument for $\gamma\Th -\delta\Ts$ shows that the top left or bottom right side
of~$\mathcal{P}_{\alpha\beta\gamma\delta}(T)$ is tangent to~$\wtclos$.

Finally, if $T= c \, Q$ for $c \in \mathbb{C}$ and $Q \in \mathcal{B}(H)$ with $Q=Q^{*}$, then the proof of the statement in the theorem
follows from the fact that $W(c\,Q)=c\,W(Q)$ and from $\sup_{\|u\|=1} | \uQu | \leq \|Q\|$, with equality for the spectral norm.
This concludes the proof.
\end{proof}

%
%

Theorem~\ref{theorem_rectangles} with an appropriate choice of the parameters $\alpha,\beta,\gamma,\delta$ 
implies the following corollary,
which leads to a polygon that contains the numerical range and exhibits useful properties.
%
%
%
%
\begin{corollary}   
\label{corollary_octagon}    
Let $T \in \mathcal{B}(H)$ have the Cartesian decomposition $T=\Th + i\Ts$, and let $T \neq c\, Q$ 
for any $c \in \mathbb{C}$ and $Q \in \mathcal{B}(H)$ with $Q=Q^{*}$.
Then $\wtclos$ is contained in a convex polygon, defined, for any norm, by the eight (not necessarily distinct) vertices
in the complex plane
\begin{eqnarray*}
& \Bigl ( \|\Ts\| - \|Th - \Ts\| \, , \, \|\Ts\| \Bigr ) \, , & \Bigl ( \|Th - \Ts\| - \|\Ts\| \, , \, -\|\Ts\| \Bigr ) \, , \\
& \Bigl ( \|Th + \Ts\|-\|\Ts\| \, , \, \|\Ts\| \Bigr ) \, , &  \Bigl ( \|\Ts\| - \|Th + \Ts\| \, , \, -\|\Ts\| \Bigr ) \; , \\
& \Bigl ( \|\Th\| \, , \, \|Th + \Ts\|-\|\Th\| \Bigr ) \, , & \Bigl ( -\|\Th\| \, , \, \|\Th\| - \|Th + \Ts\| \Bigr ) \; , \\
& \Bigl ( \|\Th\| \, , \, \|\Th\| - \|\Th - \Ts\| \Bigr ) \, , & \Bigl ( -\|\Th\| \, , \, \|\Th - \Ts\| - \|\Th\| \Bigr ) \; ,
\end{eqnarray*}
resulting in a quadrilateral, hexagon, or octagon that is symmetric with respect to the origin. If the norm is the spectral norm, then this
polygon is tangent to~$\wtclos$ in at least four points: one for each pair of opposing sides.

The numerical radius $w(T)$ of $T$ satisfies the inequality
\beq
\nonumber
w(T) \leq \Biggl ( \max \Bigl \{ \eta_{1}^{2} + \|\Th\|^{2} , \eta_{2}^{2} + \|\Ts\|^{2} \Bigr \} \Biggr )^{1/2} \; ,
\eeq
where
\begin{eqnarray*}
& & \eta_{1} = \max \Biggl \{ \Bigl | \|\Th + \Ts \| - \|\Th\| \Bigr | , \Bigl | \|\Th - \Ts \| - \|\Th\| \Bigr | \Biggr \} \; , \\ 
& & \eta_{2} = \max \Biggl \{ \Bigl | \|\Th + \Ts \| - \|\Ts\| \Bigr | , \Bigl | \|\Th - \Ts \| - \|\Ts\| \Bigr | \Biggr \}  \; .
\end{eqnarray*}
If $T=c\, Q$ for $c \in \mathbb{C}$ and $Q \in \mathcal{B}(H)$ with $Q=Q^{*}$,
then~$\wtclos$ is contained in the closed line segment determined by the endpoints $\pm c\, \|Q\|_{\sigma}$, at least one of which is 
a boundary point of~$\wtclos$ if the norm is the spectral norm. 
\end{corollary}   
\begin{proof}
Theorem~\ref{theorem_rectangles} with $\alpha=\beta=\gamma=\delta=1$ implies that~$\wtclos$ is contained 
in the intersection of the rectangle $[-\|\Th\|,\|\Th\|] \times [-\|\Ts\|,\|\Ts\|]$ and the parallelogram~$\mathcal{P}_{1111}(T)$, defined, 
after the usual identification of $\mathbb{C}$ with $\mathbb{R}^{2}$, by the lines $x \pm y  = \pm \|\Th \pm \Ts \|$,
which is tangent to it in at least four points.

Theorem~\ref{theorem_rectangles} shows that the corners of the rectangle~$\mathcal{R}(T)$ are cut off by these lines. 
If we label the top and right-hand sides of~$\mathcal{R}(T)$, respectively, as $S_{1}$ and $S_{2}$, then the vertices 
of~$\mathcal{R}(T) \cap \mathcal{P}_{1111}(T)$ are given by the following four intersection points and their reflections with respect to the origin: 
\begin{eqnarray*}
& & L_{1} \cap S_{1} = \Bigl ( \|Th + \Ts\|-\|\Ts\| \, , \, \|\Ts\| \Bigr ) \; , \\
& & L_{4} \cap S_{1} = \Bigl ( \|\Ts\| - \|Th - \Ts\| \, , \, \|\Ts\| \Bigr ) \; , \\
& & L_{1} \cap S_{2} = \Bigl ( \|\Th\| \, , \, \|Th + \Ts\|-\|\Th\| \Bigr ) \; , \\
& & L_{3} \cap S_{2} = \Bigl ( \|\Th\| \, , \, \|\Th\| - \|\Th - \Ts\| \Bigr ) \; ,
\end{eqnarray*}
where the lines $L_{j}$ are the same lines as in the proof of Theorem~\ref{theorem_rectangles} with $\alpha=\beta=\gamma=\delta=1$.
These are precisely the vertices in the statement of the corollary.
Moreover, each side of~$\mathcal{R}(T)$ contains two vertices of the intersection, which may coincide. To show this, it is sufficient
to consider $S_{1}$, as the arguments for the other sides are analogous.
The real part of the intersection of $S_{1}$ with $L_{4}$ satisfies 
\begin{multline}
\nonumber
\|\Ts\| - \|\Th - \Ts \| = \|\Ts\| - \|\Th + \Ts - 2\Ts \| \\
\leq \|\Ts\| - \lb 2\|\Ts\| - \|\Th + \Ts \| \rb = \|\Th + \Ts \| - \|\Ts\| \; ,
\end{multline}
which means that this vertex lies to the left of the intersection of $S_{1}$ with $L_{1}$, although it may coincide with it. As a result,
the intersection of $\mathcal{R}(T)$ and $\mathcal{P}_{1111}(T)$ takes the form of a quadrilateral, hexagon, or octagon.

The polygon determined by these vertices is closed and convex, since it is the intersection of two closed convex sets, so that the largest distance
from the origin to any point in the polygon is obtained at one or more vertices. Defining,
\begin{eqnarray*}
& & \eta_{1} = \max \Biggl \{ \Bigl | \|\Th + \Ts \| - \|\Th\| \Bigr | , \Bigl | \|\Th - \Ts \| - \|\Th\| \Bigr | \Biggr \} \; , \\ 
& & \eta_{2} = \max \Biggl \{ \Bigl | \|\Th + \Ts \| - \|\Ts\| \Bigr | , \Bigl | \|\Th - \Ts \| - \|\Ts\| \Bigr | \Biggr \}  \; .
\end{eqnarray*}
that maximum distance is given by
\bdis
\Biggl ( \max \Bigl \{ \eta_{1}^{2} + \|\Th\|^{2} , \eta_{2}^{2} + \|\Ts\|^{2} \Bigr \} \Biggr )^{1/2} \; ,
\edis
which is necessarily an upper bound on the numerical radius $w(T)$.

Finally, the statement in the corollary for the case $T= c \, Q$, $c \in \mathbb{C}$, and $Q \in \mathcal{B}(H)$ with $Q=Q^{*}$, follows immediately
from the corresponding case in Theorem~\ref{theorem_rectangles}.
This concludes the proof.
\end{proof}

Figure~\label{fig2} illustrates Corollary~\ref{corollary_octagon} for the matrices
\small
\bdis
A=
\begin{pmatrix}
2-4i & -4+4i & -4-i  & 3-i   \\
1-3i & -1    & -2+2i & 5i    \\
- 2i & 4-i   & -1-2i & 3-4i  \\
4-4i & -4i   &  1-3i & 2+5i  \\
\end{pmatrix}
\;\; \text{and} \;\;
B=
\begin{pmatrix}
4- i  & -3+2i & 3+5i  & -2+3i \\
-1-i  & 0     & 1-4i  & -3-2i \\
-4+4i & 1-4i  & -4i   & -2    \\
4+i   & 4+i   & 2+2i  & 1     \\
\end{pmatrix}
\; ,
\edis
\normalsize
where the solid outer circle represents the bound from~(\ref{bnd_classical}), the solid inner circle shows the bound from~(\ref{kittaneh_ineq2}), and the dashed
circle represents the bound from~(\ref{kittaneh_ineq1}). The polygon is the one obtained from Corollary~\ref{corollary_octagon} and the shaded area is 
the numerical range of the matrix. The octagon can clearly either be a very good approximation to the \textit{numerical range} as for the matrix~$A$ or 
it can be less satisfactory as for the matrix~$B$. However, in both cases, the approximation to the \textit{numerical radius} is equally good. The latter 
remains true even for very elongated numerical ranges with an area much smaller than that of the approximating octagon.
%
%
%
%
\begin{figure}[H]
\begin{center}
\raisebox{0ex}{\includegraphics[width=0.38\linewidth]{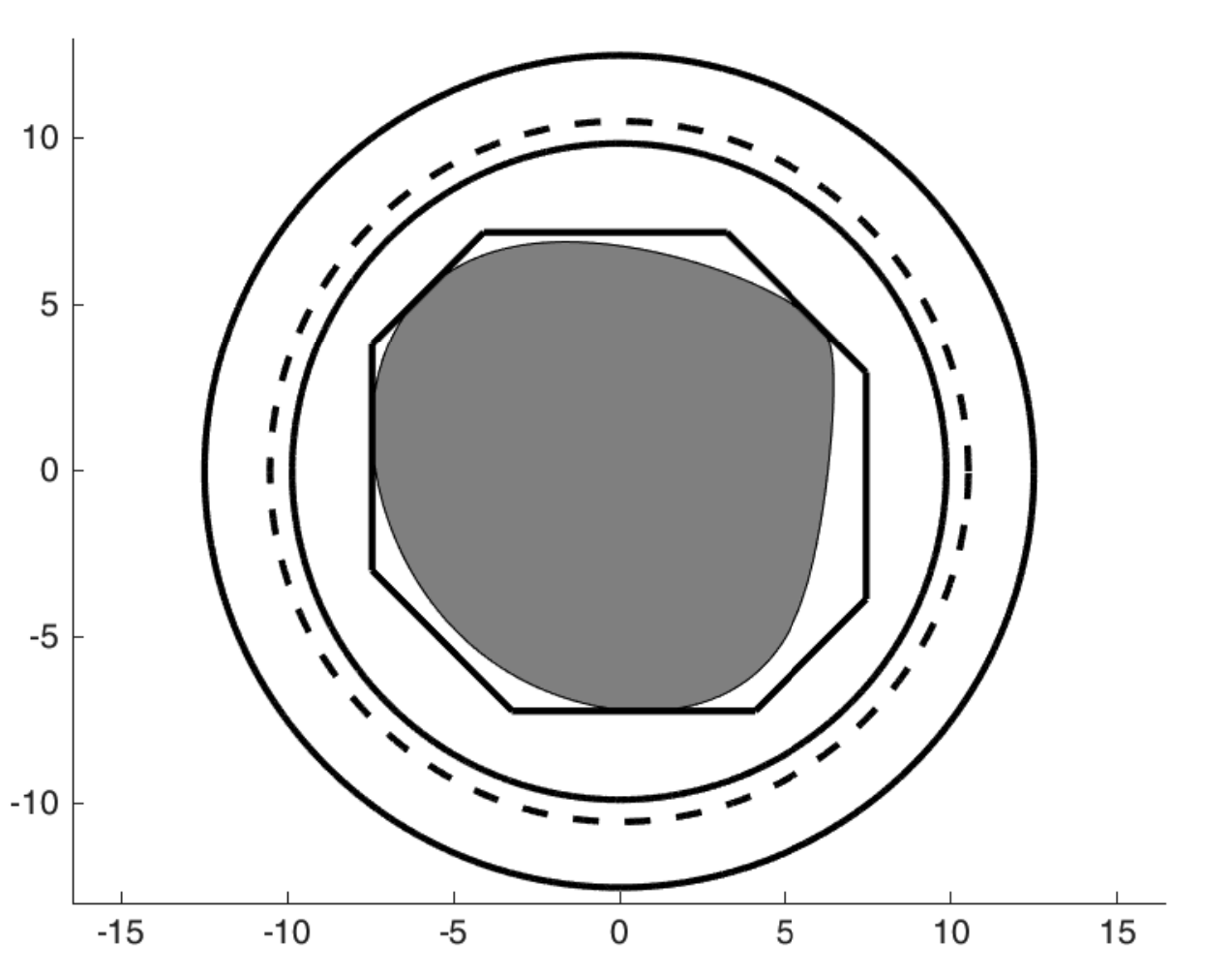}}
\hskip 1.00cm  
\raisebox{-1.25ex}{\includegraphics[width=0.38\linewidth]{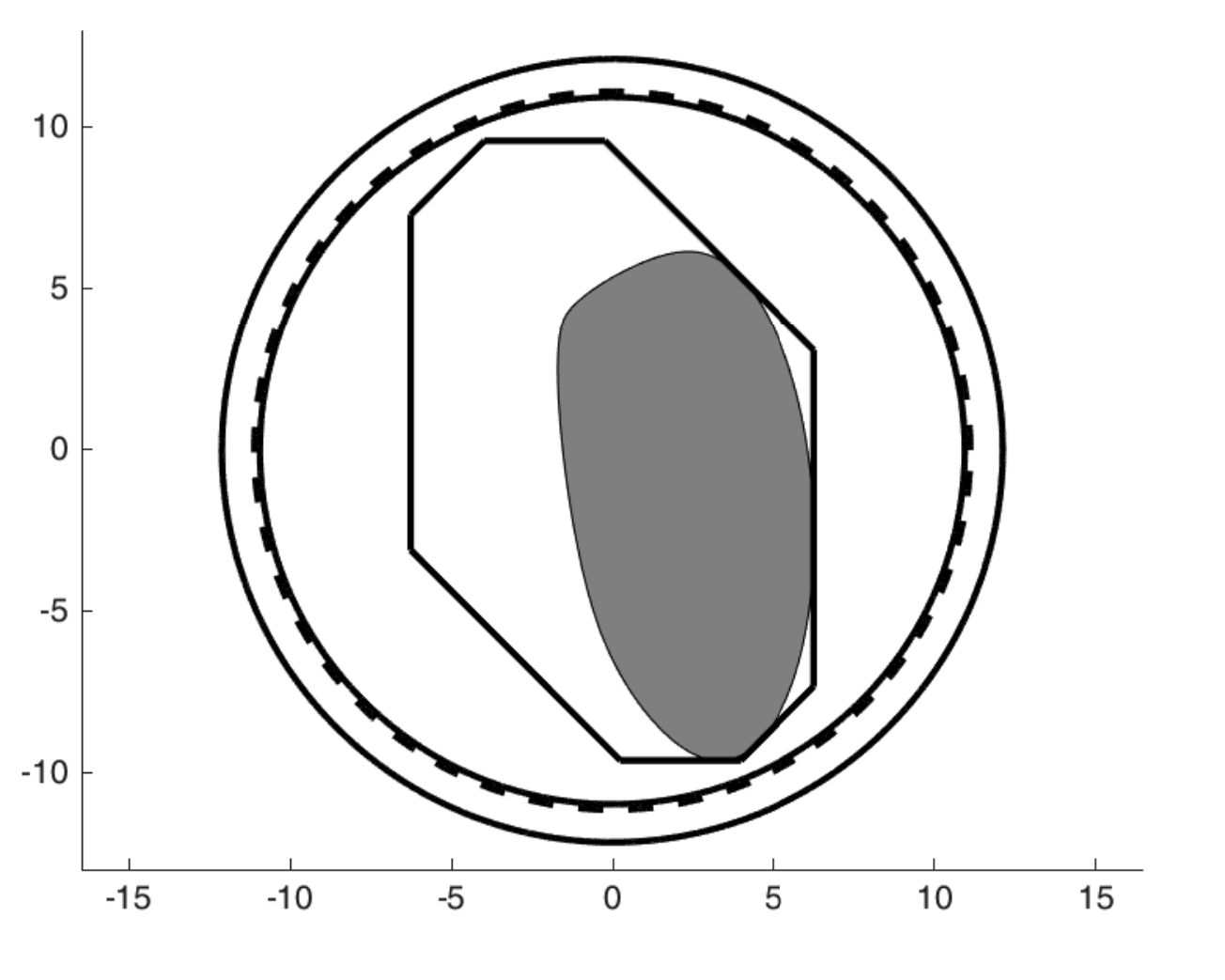}}
\caption{Octagons containing the numerical ranges of the matrices~$A$ and~$B$.}
\label{fig2}                      
\end{center}
\end{figure}

The bound on the numerical radius obtained in Corollary~\ref{corollary_octagon} is not necessarily better than existing bounds, although it often is.
To obtain an idea of the relative performance of the bound in Corollary~\ref{corollary_octagon} with the \emph{spectral norm}, we have compared it to the 
bounds in~(\ref{kittaneh_ineq1}) and~(\ref{kittaneh_ineq2}). 
To do this, we have generated $1000$ $m \times m$ matrices, with $m=10, 100, 500, 1000$, whose elelements are complex with real and complex 
parts uniformly randomly distributed in the interval $[-4,4]$. We have listed in Table~\ref{table1}, the average ratios of the respective bounds
to the spectral norm of the matrix (the smaller the ratio, the better the bound), which demonstrates the advantage of Corollary~\ref{corollary_octagon}. 
Moreover, the results appear to be quite insensitive to the size of the matrix.

%
%
\begin{table}[H]
\begin{center}
\small           
\begin{tabular}{c|c c c}
m      &  Bound~(\ref{kittaneh_ineq1})     &  Bound~(\ref{kittaneh_ineq2})  &  Corollary~\ref{corollary_octagon}      \\ \hline 
       &                                 &                              &                                       \\
10     &  0.91                           &  0.88                        &  0.80                                 \\
100    &  0.90                           &  0.86                        &  0.77                                 \\
500    &  0.90                           &  0.86                        &  0.77                                 \\
1000   &  0.90                           &  0.85                        &  0.77                                 \\
\end{tabular}
\caption{Comparison of bounds on the numerical range for $m=10,100,500,1000$.}
\label{table1}
\end{center}
\end{table}

\end{document}